\documentclass{article}
\usepackage{graphicx} % Required for inserting images
\usepackage[margin=1in]{geometry}
\usepackage[utf8]{inputenc}
\usepackage{amsmath}
\usepackage{amsthm}
\usepackage{amsfonts}
\usepackage{amssymb}
\usepackage{ragged2e}
\usepackage{authblk}
\usepackage{cite}
\usepackage{enumerate}
\usepackage{hyperref}

\providecommand{\keywords}[1]
{
  \small	
  \textbf{\textit{Keywords: }} #1
}

\newtheorem{theorem}{Theorem}
\newtheorem{corollary}[theorem]{Corollary}

\newtheorem{proposition}[theorem]{Proposition}

\title{Revisiting an infinitely nested radical}
\author{Aung Phone Maw}
\date{February 2026}

\begin{document}

\maketitle

\begin{abstract}
    We revisit an infinitely nested radical by Ramanujan. Utilizing the full strength of his method, we shall arrive at some new infinitely nested radicals.\\
\end{abstract}

\section{Introduction}

Ramanujan posed the following problem in the Journal of Indian Mathematical Society in 1911, the problem (\cite[p.~323]{Ramanujan}) is to find the value of :

\[
\sqrt{1+2\sqrt{1+3\sqrt{1+4\sqrt{1+\cdots}}}}. \tag{1.1} \label{exp1}
\]
The value of which is simply $3$. The following algebraic relation is used by Ramanujan to arrive at the solution of the above infinitely nested radical : 

\[
x^2= n^2+an+ax + (x-n-a)(x+n). \tag{1.2} \label{ramanujanrelation}
\]
Ramanujan used this relation to obtain a general infinitely nested radical and then obtaining the value of \eqref{exp1} as a particular case (see \cite[p.~108]{Berndt}).
Utilizing the full strength of his method we shall evaluate the following infinitely nested radicals :

\[
3=\sqrt{1^2+(5+1^2-4)\sqrt{3^2+(5+2^2-8)\sqrt{5^2+(5+3^2-12)\sqrt{7^2+\cdots}}}} \tag{1.3}
\]

\[
3=\sqrt{2^2+(3+1^2-3)\sqrt{4^2+(3+2^2-6)\sqrt{6^2+(3+3^2-9)\sqrt{8^2+\cdots}}}} \tag{1.4}
\]
\[
9=\sqrt[3/2]{8+3+1\ \sqrt[3/2]{8+6+2\ \sqrt[3/2]{8+9+3\ \sqrt[3/2]{8+12+4\ \sqrt[3/2]{8+15+\cdots}}}}} \tag{1.5}
\]
\[
1 = \sqrt{\sqrt{2}-1+(\sqrt{2}-1)\sqrt{\sqrt{3}-1+(\sqrt{3}-1)\sqrt{\sqrt{4}-1+(\sqrt{4}-1)\sqrt{\cdots}}}}. \tag{1.6}
\]

\section{Three Identities}
Let $a_1,a_2,...$ and $n_1,n_2,..$ be two sequences, then by \eqref{ramanujanrelation}, for all $x$ and for each $i \in \mathbb{N}$ we have :
\[
x^2=n_i^2+a_in_i+a_ix+(x-n_i-a_i)(x+n_i). \tag{2.1}
\]
Thus for all $x$ and for all $i\in \mathbb{N}$ we have :
\[
F(x)^2= n_i^2+a_in_i+a_ix+(x-n_i-a_i)F(x+n_i),
\]
with $F(x)=x$. Therefore, we have the following formal iteration :
\begin{flalign*}
    &F(x) \\&= \sqrt{n_1^2+a_1n_1+a_1x+(x-n_1-a_1)F(x+n_1)}\\
    &= \sqrt{n_1^2+a_1n_1+a_1x+(x-n_1-a_1)\sqrt{n_2^2+a_2n_2+a_2(x+n_1)+(x+n_1-n_2-a_2)F(x+n_2)}}\\
    & \vdots\\
    &=\sqrt{n_1^2+a_1n_1+a_1x+(x-n_1-a_1)\sqrt{n_2^2+a_2n_2+a_2(x+n_1)+(x+n_1-n_2-a_2)\sqrt{ n_3^2+a_3n_3 + \cdots }}}
\end{flalign*}

Thus :

\begin{proposition}
    Let $a_1,a_2,...$ and $n_1,n_2,..$ be two sequences, formally, we have :
    \[
    x = \sqrt{n_1^2+a_1n_1+a_1x+(x-n_1-a_1)\sqrt{n_2^2+a_2n_2+a_2(x+n_1)+(x+n_1-n_2-a_2)\sqrt{ n_3^2+a_3n_3 + \cdots }}},
    \tag{2.2} \label{identity1}\]
    for all $x$.
\end{proposition}

Another identity for degree three can be obtained by the following algebraic relation:
\[
x^3= 2n^3+3n^2x+(x-2n)(x+n)^2 \tag{2.3} \label{relation1}
\]
Thus for a sequence $n_1,n_2,...$ and $G(x)=x^2$ we have the following formal iteration :

\begin{flalign*}
    &G(x) \\
    &= \sqrt[3/2]{2n_1^3+3n_1^2x+(x-2n_1)G(x+n_1)}\\
    &= \sqrt[3/2]{2n_1^3+3n_1^2x+(x-2n_1) \sqrt[3/2]{2n_2^3+3n_2^2(x+n_1)+(x+n_1-2n_2)G(x+n_2)}}\\
    &\vdots \\
    &= \sqrt[3/2]{2n_1^3+3n_1^2x+(x-2n_1) \sqrt[3/2]{2n_2^3+3n_2^2(x+n_1)+(x+n_1-2n_2)\sqrt[3/2]{2n_3^3+3n_3^2(x+n_1+n_2)+\cdots}}}
\end{flalign*}
Therefore :
\begin{proposition} Let $n_1,n_2,...$ be a sequence, formally, we have :
\[
x^2= \sqrt[3/2]{2n_1^3+3n_1^2x+(x-2n_1) \sqrt[3/2]{2n_2^3+3n_2^2(x+n_1)+(x+n_1-2n_2)\sqrt[3/2]{2n_3^3+3n_3^2(x+n_1+n_2)+\cdots}}}, \tag{2.4} \label{identity2}
\]
for all $x$.
    
\end{proposition}

Finally, we consider the following relation :

\[
x^2 = a\sqrt{x^2+n}-n+(\sqrt{x^2+n}-a)\sqrt{x^2+n}, \tag{2.5} \label{relation2}
\]
Thus, with $H(x)=x$ and sequences $a_1,a_2,...$ and $n_1,n_2,...$, formally we have :
\begin{flalign*}
    &H(x)\\
    &= \sqrt{a_1\sqrt{x^2+n_1}-n_1+(\sqrt{x^2+n_1}-a_1)H(\sqrt{x^2+n_1})}\\
    &= \sqrt{a_1\sqrt{x^2+n_1}-n_1+(\sqrt{x^2+n_1}-a_1)\sqrt{a_2\sqrt{x^2+n_1+n_2}-n_2+(\sqrt{x^2+n_1+n_2}-a_2)H(\sqrt{x^2+n_1+n_2})}}\\
    &\vdots\\
    &=\sqrt{a_1\sqrt{x^2+n_1}-n_1+(\sqrt{x^2+n_1}-a_1)\sqrt{a_2\sqrt{x^2+n_1+n_2}-n_2+(\sqrt{x^2+n_1+n_2}-a_2)\sqrt{\cdots}}}.
\end{flalign*}
Therefore :
\begin{proposition}Let $a_1,a_2,...$ and $n_1,n_2,..$ be two sequences, formally, we have:
\[
x = \sqrt{a_1\sqrt{x^2+n_1}-n_1+(\sqrt{x^2+n_1}-a_1)\sqrt{a_2\sqrt{x^2+n_1+n_2}-n_2+(\sqrt{x^2+n_1+n_2}-a_2)\sqrt{\cdots}}},
\tag{2.6} \label{identity3}\]
for all $x$.    
\end{proposition}

Now we turn to particular cases. Put $x=3,a_1=a_2=\cdots=0$ and $n_i = 2i-1,$ for all $i\in \mathbb{N}$ in \textbf{Proposition 1}, then since $3+(\sum_{i=1}^{n-1}(2i-1))-(2n-1)= 3+(n-1)^2-2n+1 = 5+n^2-4n$, for each $n \in \mathbb{N}$, we have :
\begin{corollary} 
    \[
    3=\sqrt{1^2+(5+1^2-4)\sqrt{3^2+(5+2^2-8)\sqrt{\cdots\sqrt{(2n-1)^2+(5+n^2-4n)\sqrt{\cdots}}}}}. \]
\end{corollary}
Next, put $x=3,a_1=a_2=\cdots=0$ and $n_i = 2i$ for all $i \in \mathbb{N}$ in \textbf{Proposition 1}. Since $3+\sum_{i=1}^{n-1}2i - 2n = 3+n^2-n-2n = n^2-3(n-1)$, we have :

\begin{corollary}
    \[
3=\sqrt{2^2+1^2\sqrt{4^2+(2^2-3)\sqrt{\cdots\sqrt{(2n)^2+(n^2-3(n-1))\sqrt{\cdots}}}}}.
\]
\end{corollary}
In \textbf{Proposition 2}, let $x=3, n_1=n_2=\cdots=1,$ to get 
\begin{corollary}
    \[
9=\sqrt[3/2]{8+3+1\ \sqrt[3/2]{8+6+2\ \sqrt[3/2]{\cdots\sqrt[3/2]{8+3n+n\ \sqrt[3/2]{\cdots}}}}}.
\]
\end{corollary}
Finally, let $x=1, a_1=a_2=\cdots=1$ and $n_1=n_2=...=1$ in \textbf{Proposition 3} to get:
\begin{corollary}
    \[
1 = \sqrt{\sqrt{2}-1+(\sqrt{2}-1)\sqrt{\sqrt{3}-1+(\sqrt{3}-1)\sqrt{\sqrt{4}-1+(\sqrt{4}-1)\sqrt{\cdots}}}}. 
\]
\end{corollary}

\end{document}